\documentclass[11pt,fleqn]{article}
\usepackage{amssymb,latexsym,amsmath,amsfonts}
\usepackage{graphicx}

    \advance\textheight by \topskip

    \numberwithin{equation}{section}

    \DeclareMathOperator*{\IM}{Im}
    
    \def\Im{\IM}

    \newtheorem{theorem}{Theorem}[section]
    \newtheorem{lemma}[theorem]{Lemma}
    \newtheorem{corollary}[theorem]{Corollary}

    \newtheorem{Definition}[theorem]{Definition}
    \newenvironment{definition}{\begin{Definition}\rm}{\end{Definition}}
    \newtheorem{Remark}[theorem]{Remark}
    \newenvironment{remark}{\begin{Remark}\rm}{\end{Remark}}
    \newtheorem{Example}[theorem]{Example}

    \newenvironment{proof}%
    {\rm \trivlist \item[\hskip \labelsep{\bf Proof. }]}%
    {\hspace*{\fill}$\Box$\endtrivlist}
    {\rm \trivlist \item[\hskip \labelsep{\bf Proof}]}%
    {\hspace*{\fill}$\Box$\endtrivlist}

\begin{document}
 \begin{center} \Large\bf
 Multiple orthogonal polynomials of mixed type  and
 non-intersecting Brownian motions
\end{center}

    \begin{center}  \large
        E. Daems \\
            \normalsize \em
            Department of Mathematics, Katholieke Universiteit Leuven, \\
            Celestijnenlaan 200 B, 3001 Leuven, Belgium \\
            \rm evi.daems@wis.kuleuven.be \\[3ex]
            \rm and \\[3ex]
            \large
        A.B.J. Kuijlaars \\
            \normalsize \em
            Department of Mathematics, Katholieke Universiteit Leuven, \\
            Celestijnenlaan 200 B,  3001 Leuven, Belgium \\
            \rm arno@wis.kuleuven.be
     \end{center}\ \\[1ex]

\begin{abstract}
We present a generalization of multiple orthogonal polynomials of
type I and type II, which we call multiple orthogonal polynomials
of mixed type. Some basic properties are formulated, and a Riemann-Hilbert
problem for the multiple orthogonal polynomials of mixed type is given. We
derive a Christoffel-Darboux formula for these polynomials using
the solution of the Riemann-Hilbert problem.  The main motivation
for studying these polynomials comes from a model of non-intersecting
one-dimensional Brownian motions with a given number of starting points
and endpoints. The correlation kernel for the positions of the Brownian paths
at any intermediate time coincides with the Christoffel-Darboux kernel
for the multiple orthogonal polynomials of mixed type with respect
to Gaussian weights.
\end{abstract}

\section{Introduction}
In the early nineties, Fokas, Its and Kitaev \cite{FIK} introduced
a $2 \times 2$ matrix valued Riemann-Hilbert problem that characterizes
orthogonal polynomials on the real line. This approach can be used to study
various aspects of the theory of orthogonal polynomials.
Combined with the Deift-Zhou steepest descent method for Riemann-Hilbert
problems it has been very succesful in deriving asymptotics for orthogonal polynomials
and solving basic questions in the theory of random matrices \cite{Dei}.
The Deift-Zhou steepest descent method was first introduced in
\cite{Dei2} and further developed in for example \cite{DKMVZ1} and
\cite{DKMVZ2}.

In \cite{VAGK}, this Riemann-Hilbert problem was extended to the
case of multiple orthogonal polynomials of type I and type II.
Multiple orthogonal polynomials are polynomials which satisfy
orthogonality conditions with respect to a number of measures. The
definition will be given in section 2. The Riemann-Hilbert problem
is now of size $(p+1)\times (p+1)$, where $p$ is the number of measures.

In this paper, we introduce a generalization of multiple orthogonal
polynomials. These new polynomials satisfy orthogonality
conditions with respect to two sets of weights
$w_{1,1},\ldots,w_{1,p}$ and $w_{2,1},\ldots,w_{2,q}$.
We call these polynomials  multiple orthogonal polynomials of mixed type. The
definition will be given in section 2, together with some basic
properties concerning existence and uniqueness. We can again
characterize these polynomials by a Riemann-Hilbert problem which is
now of size $(p+q)\times (p+q)$. This will be given in section 3.

The usual monic orthogonal polynomials $P_n$ on the real line with
weight function $w$ satisfy a three term recurrence relation
and this gives rise to the basic Christoffel-Darboux formula (see
for example \cite{Chih})
\begin{align}\label{christoffel0}
\sum_{j=0}^{n-1}\frac{1}{h_j}P_j(x)P_j(y) =
    \frac{1}{h_{n-1}}\frac{P_{n}(x)P_{n-1}(y)-P_{n-1}(x)P_{n}(y)}{x-y},
\end{align}
where
\begin{align}
 h_j = \int P_j(x) x^j w(x)dx.
\end{align}
This formula was generalized to multiple orthogonal polynomials in
\cite{BK1} in the case of two weights and in \cite{DK} in the
general case of $p$ weights. We will derive using the
Riemann-Hilbert problem a Christoffel-Darboux formula for the
multiple orthogonal polynomials of mixed type.

Our main motivation for studying this new kind of orthogonality
comes from non-intersecting Brownian paths. Consider $n$
independent one-dimen\-sional Brownian motions that start in $n$
different fixed points at time $t=0$ and end in $n$ different
fixed points at time $t=1$, conditioned on the fact that they do
not intersect in the full time interval $(0,1)$. At any
intermediate time $t \in (0,1)$ the positions of the Brownian
paths are distributed according to a determinantal point process on the
real line.
This is a consequence of a classical theorem of Karlin and McGregor
\cite{KMcG} and it applies not only to Brownian motion, but to any
one-dimensional strong Markov process with continuous sample
paths.

Of special interest is the confluent case in which
many of the starting points and many of the endpoints coincide.
The formulas simplify for Brownian motion because of
the properties of the Gaussian transition probabilities.
In the extreme case of one starting point and one endpoint,
the positions of the Brownian paths have the
same distribution, up to simple scaling, as the eigenvalues of a
matrix from the Gaussian Unitary Ensemble (GUE) which is a basic ensemble
from random matrix theory. In this case the kernel for the
determinantal point process is constructed out of Hermite polynomials
and the Christoffel-Darboux formula (\ref{christoffel0}) expresses
this kernel in terms of Hermite polynomials of degrees $n$ and $n-1$ only.

In the case of one starting point and $q$ endpoints the positions
of the Brownian paths have the same distribution as the
eigenvalues of a Gaussian unitary matrix with external source \cite{ABK3,AvM,TW}.
Then the kernel is constructed out of multiple Hermite polynomials
of type I and II. The Christoffel-Darboux formula for multiple
orthogonal polynomials of \cite{DK} expresses this kernel in terms
of a sum of $q+1$ terms, in which each term involves products of
multiple Hermite polynomials of type I and type II.

The next step is to consider Brownian motions with begin in $p$
starting points and end in $q$ endpoints. This gives rise to the
multiple orthogonal polynomials of mixed type with respect
to Gaussian weights. We will discuss this in more detail in
section 6. Unfortunately we do not know if there
exists a corresponding random matrix model

For a summary of the above discussion, see Table \ref{tabel1}.

\newpage

\begin{table}\label{tabel1}
    \begin{center}
    \begin{tabular}{|p{3.5cm} |p{3.5cm}| p{3.5cm}|}\hline
         Non-intersecting Brownian motions  &   Associated polynomials &
         Random matrix ensemble
        \\ \hline \hline
        \begin{center}1 starting point and 1 endpoint  \end{center} &
        \begin{center} Hermite polynomials \end{center} &
        \begin{center} Gaussian unitary ensemble
        \end{center}
        \\\hline
        \begin{center}1 starting point and $q \geq 2$ endpoints \end{center}&
        \begin{center} multiple Hermite polynomials \end{center} &
        \begin{center} Gaussian unitary ensemble  with external source\end{center}
        \\\hline
        \begin{center} $p \geq 2$ starting points and $q \geq 2$ endpoints  \end{center} &
        \begin{center}multiple Hermite polynomials of mixed type \end{center} &
        \begin{center} unknown \end{center}
        \\\hline
    \end{tabular}
    \end{center}
    \caption{Overview of the connection between non-intersecting Brownian motions, associated polynomials,
    and random matrix ensembles. The three cases are illustrated in Figures \ref{figuur1},
    \ref{figuur2}, and \ref{figuur3}. See section 6 for a more
    detailed discussion.}
\end{table}

\newpage

\begin{figure}[ht]
    \begin{center}
        \includegraphics[width=8cm,height=6cm]{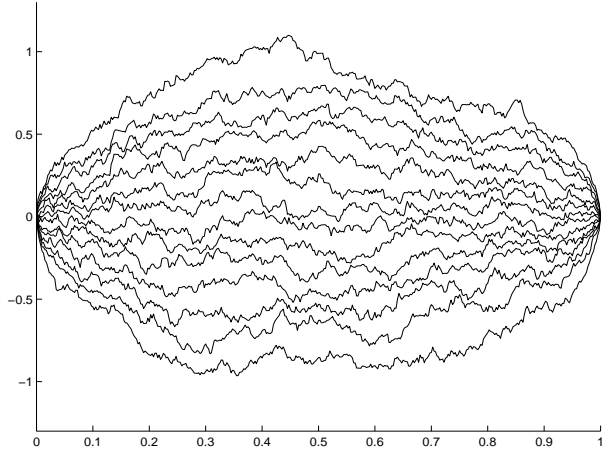}
        \caption{Non intersecting Brownian motions which start and end at 1 point.
        At any intermediate time the positions of the paths have the same distribution
        as the eigenvalues of a GUE matrix. The correlation kernel is built out
        of Hermite polynomials.}\label{figuur1}
    \end{center}
\end{figure}

\begin{figure}[ht]
    \begin{center}
        \includegraphics[width=8cm,height=6cm]{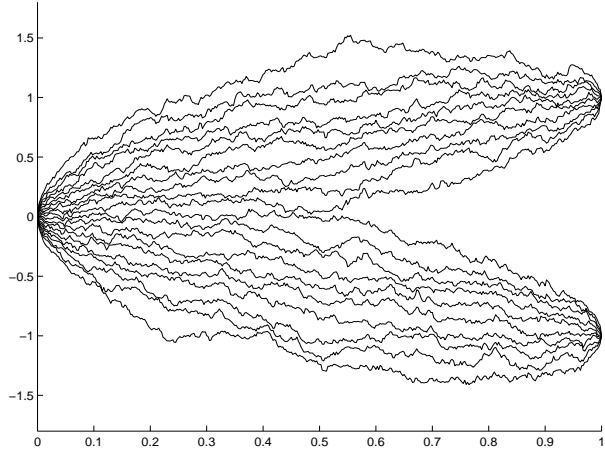}
        \caption{Non intersecting Brownian motions which start at 1 point and end at 2 different points.
        At any intermediate time the positions of the paths have the same distribution
        as the eigenvalues of a Gaussian unitary random matrix with external source.
        The correlation kernel is built out of multiple Hermite polynomials.}\label{figuur2}
    \end{center}
\end{figure}

\begin{figure}[ht]
    \begin{center}
        \includegraphics[width=8cm,height=6cm]{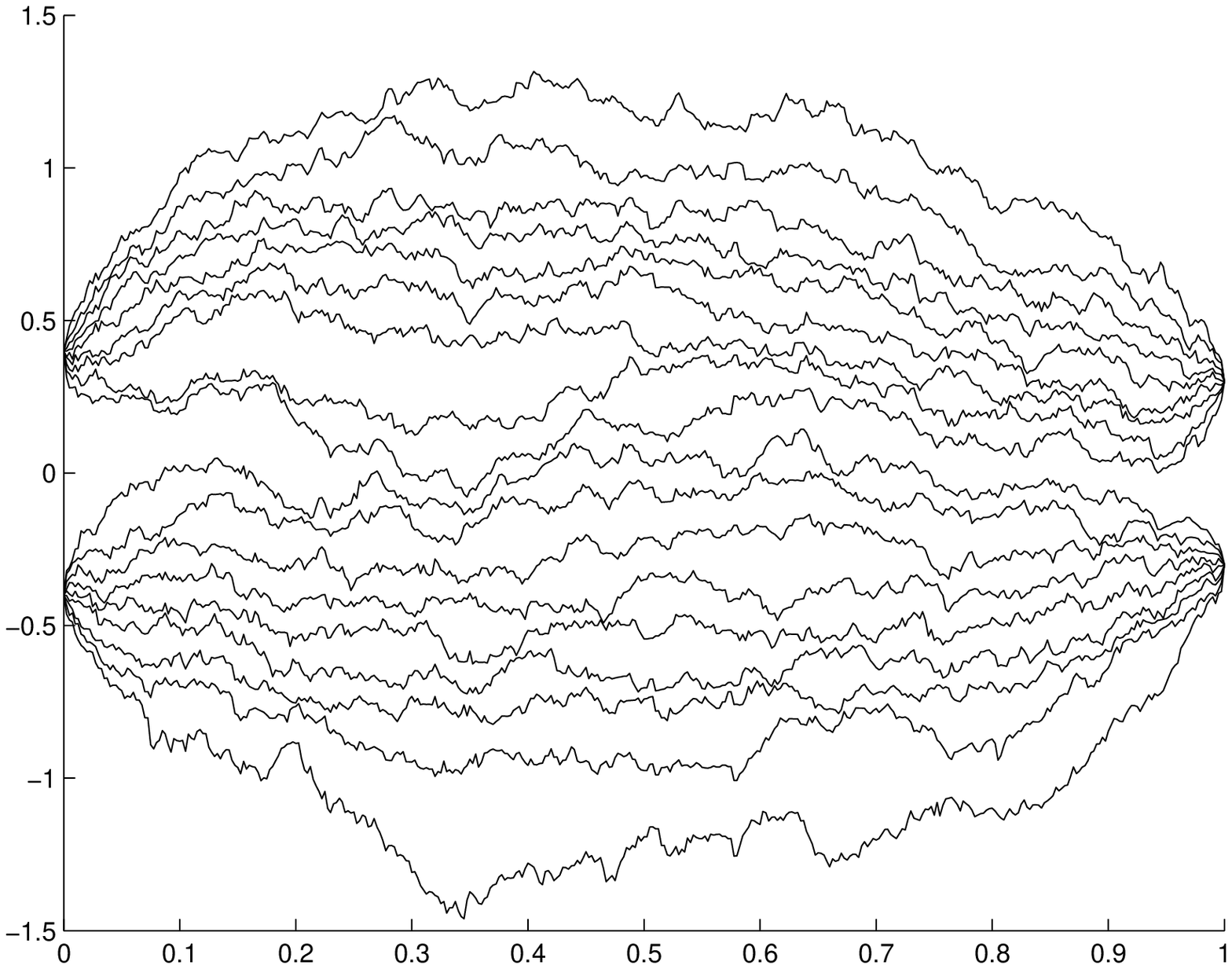}
        \caption{Non intersecting Brownian motions which start at 2 points and end at 2 points.
        At any intermediate time the positions of the paths are distributed according to
        a determinantal point process with a kernel that is built out of
        multiple Hermite polynomials of mixed type.}\label{figuur3}
    \end{center}
\end{figure}

\section{Multiple orthogonal polynomials of mixed type}
In this section we define the multiple orthogonal
polynomials of mixed type and we give conditions for existence.
These polynomials can be seen as a generalization of
multiple orthogonal polynomials of type I and type II which we discuss
first. Throughout this paper we will say that
$w$ is a weight on $\mathbb R$ if $w(x) \geq 0$ for $x \in \mathbb R$
and  $\int x^k w(x) dx < +\infty$ for every $k \in \mathbb N \cup \{0\}$.

\subsection{Multiple orthogonal polynomials of type I}

Let $w_1, w_2, \ldots, w_p$ be $p$ weights on the real line and let
$\vec{n} = (n_1,\ldots,n_p)$ be a multi-index consisting of non-negative
integers. If $A_1, \ldots, A_p$ are polynomials and
\begin{align} \label{linform}
    Q(x) = \sum_{j=1}^p A_j(x) w_j(x), \qquad \mbox{deg } A_j \leq n_j-1,
\end{align}
such that
\begin{align} \label{typeIdef}
    \int Q(x) x^j dx  = 0 \qquad \mbox{for } j=0, \ldots, |\vec{n}|-2,
\end{align}
then the $A_j$ are called multiple orthogonal polynomials of type I
and $Q$ is the linear form built out of the multiple orthogonal
polynomials of type I. Here we follow the usual multi-index notation
\[ |\vec{n}| = \sum_{i=1}^p n_i. \]
The relations (\ref{typeIdef}) give us $|\vec{n}|-1$ homogeneous
linear equations for the in total $|\vec{n}|$ coefficients of the
polynomials $A_j$. So there is always a non-zero solution. If the
solution is unique up to a multiplicative factor, then the
multi-index $\vec{n}$ is called normal for type I. The multi-index
$\vec{n}$ is called strongly normal for type I if we have
\[ \int Q(x) x^{|\vec{n}|-1} dx \neq 0 \]
for any non-zero $Q$ satisfying (\ref{linform})--(\ref{typeIdef}).
In that case we can normalize the multiple orthogonal polynomials
of type I so that
\begin{align} \label{typeInorm}
    \int Q(x) x^{|\vec{n}|-1} dx = 1.
\end{align}
We call (\ref{typeInorm}) a type I normalization.

\subsection{Multiple orthogonal polynomials of type II}

Let $w_1, w_2, \ldots, w_q$ be $q$ weights on the real line and let
$\vec{m} = (m_1,\ldots,m_q)$ be a multi-index of length $q$.
If $P$ is a polynomial of degree $|\vec{m}|$ such that
\begin{align} \label{typeIIdef}
    \int P(x) x^j w_k(x) dx = 0 \qquad
    \mbox{for } j=0, \ldots, m_k-1 \mbox{ and } k = 1,\ldots, q,
    \end{align}
then $P$ is called a multiple orthogonal polynomial of type II.
The $|\vec{m}|$ equations (\ref{typeIIdef}) are homogeneous linear
equations for the $|\vec{m}| + 1$ coefficients of $P$. So there is
always a non-zero solution. If the solution is unique up to a
multiplicative factor then the multi-index $|\vec{m}|$ is called
normal for type II. If every non-zero solution has a non-zero
leading coefficient then the multi-index $|\vec{m}|$ is called
strongly normal for type II. In that case we can normalize the
multiple orthogonal polynomial of type II so that
\begin{align} \label{typeIInorm}
    P_n(x) = x^n + \cdots
    \end{align}
and we call (\ref{typeIInorm}) a type II normalization.

\addvspace{1mm}
 For more details and examples of multiple orthogonal polynomials of type
I and type II, and  about their relation with Hermite-Pad\'e
approximation, we refer the interested reader to \cite{Apt,NS,VAC} and the
references cited therein.

\subsection{Multiple orthogonal polynomials of mixed type}

To define the multiple orthogonal polynomials of mixed type, we need two
sets of weights on $\mathbb R$: $w_{1,1}, w_{1,2}, \ldots,
w_{1,p}$ and $ w_{2,1}, w_{2,2}, \ldots, w_{2,q}$, which we collect in
two row vectors
\[ \vec{w}_1 = (w_{1,1}, \ldots, w_{1,p}),
    \qquad \vec{w}_2 = (w_{2,1}, \ldots, w_{2,q}), \]
and two multi-indices $\vec{n} = (n_1, \ldots, n_p)$ and
$\vec{m} = (m_1, \ldots, m_q)$
of length $p$ and $q$, respectively, such that
\begin{align} \label{normcond}
    |\vec{n}| = |\vec{m}| + 1.
\end{align}

\begin{definition}
We call the polynomials $A_1, \ldots, A_p$  with
\begin{align} \label{typeMdeg}
    \mbox{deg } A_j \leq n_{j}-1 \qquad \mbox{for } j=1,\ldots, p
    \end{align}
multiple orthogonal polynomials of mixed type for the pair of
multi-indices $(\vec{n},\vec{m})$ and with respect to the vectors
of weights $\vec{w}_1$ and $\vec{w}_2$ if the function
\begin{align} \label{typeMform}
    Q(x) = \sum_{j=1}^p A_j(x) w_{1,j}(x)
\end{align}
satisfies the following orthogonality conditions:
\begin{align} \label{typeMdef}
    \int Q(x) x^j w_{2,k}(x) dx = 0 \quad
    \mbox{ for } j=0,1, \ldots, m_k-1 \mbox{ and } k=1,\ldots,q.
\end{align}
\end{definition}
To emphasize the dependence on the multi-indices we also write
\[ A_j = A_{j,\vec{n}, \vec{m}}, \qquad Q = Q_{\vec{n},\vec{m}}, \]
and to emphasize the role of the two vectors of weights we will
occasionally write
\[ A_j(x) = A_{j, \vec{n}, \vec{m}}(x; \vec{w}_1, \vec{w}_2),
    \qquad Q(x) = Q_{\vec{n}, \vec{m}}(x; \vec{w}_1, \vec{w}_2), \]
although mostly we drop the explicit mentioning of the weights.
Note that the role of the two vectors of weights is not symmetric.
The function $Q$ is a linear form with respect to the weights from
$\vec{w}_1$ as in multiple orthogonality of type I, and the linear
form has a number of orthogonality conditions with respect to the weights from $\vec{w}_2$
as in multiple orthogonality of type II.

The conditions (\ref{typeMdef}) lead to $|\vec{m}|$ homogeneous linear
equations for the in total $|\vec{n}|$ free coefficients of the polynomials
$A_j$. Because of the assumption (\ref{normcond}) there is always a
non-zero solution. If the polynomials $A_j$ are unique up to a multiplicative
constant, then we call $(\vec{n}, \vec{m})$ a \textbf{normal} pair of indices for
the two sets of weights $\vec{w}_1$ and $\vec{w}_2$.

For a normal pair of indices we can choose a certain normalization
in order to define a unique multiple orthogonal polynomial of mixed type.
In this paper, we are going to use two types of normalizations:
\begin{itemize}
\item \textbf{Type I normalization:} Fix $k=1,\ldots, q$ and normalize
$Q$ such that
\begin{align}
 \int Q(x) x^{m_k} w_{2,k}(x) dx  = 1.
\end{align}
If we choose this normalization, we write
\begin{align}\label{normI}
A_j = A_j^{(I,k)}, \qquad Q = Q^{(I,k)},
\end{align}
or in full notation, if we want to emphasize the dependence on the multi-indices and the weights
\[ A_j(x) = A_{j,\vec{n}, \vec{m}}^{(I,k)}(x; \vec{w}_1, \vec{w}_2),
    \qquad Q(x) = Q_{\vec{n}, \vec{m}}^{(I,k)}(x; \vec{w}_1, \vec{w}_2). \]
 \item \textbf{Type II normalization:} Fix $k=1, \ldots, p$ and normalize $Q$
 such that $A_k$ is a monic polynomial of degree $n_k-1$. If we choose this
 normalization, we write
 \begin{align}\label{normII}
 A_j = A_j^{(II,k)}, \qquad Q = Q^{(II,k)},
 \end{align}
 or in full
 \[ A_j(x) = A_{j,\vec{n}, \vec{m}}^{(II,k)}(x; \vec{w}_1, \vec{w}_2),
    \qquad Q(x) = Q_{\vec{n}, \vec{m}}^{(II,k)}(x; \vec{w}_1, \vec{w}_2). \]
\end{itemize}
We emphasize that the above normalizations may
not always be possible.

\subsection{Conditions for normality}

We will state the conditions using the Hilbert space
geometry of $L^2(\mathbb R)$ and to do so we assume that
\begin{align} \label{spanF}
    x^j w_{1,k} \in L^2(\mathbb R), \qquad
    \mbox{ for } k=1, \ldots, p, \quad j = 0, 1, \ldots, n_k-1.
\end{align}
\begin{align} \label{spanG}
    x^j w_{2,k} \in L^2(\mathbb R), \qquad
    \mbox{ for } k=1, \ldots, q, \quad j=0, 1, \ldots, m_k-1
    \end{align}
Associated with $\vec{n}$ and $\vec{w}_1$ we have the vector space
\begin{align} \label{F}
    F_{\vec{n}} = \{\sum_{j=1}^p A_j w_{1,j} \mid
    A_j \mbox{ is a polynomial of degree } \leq n_j-1 \mbox{ for } j=1, \ldots, p \}
\end{align}
and associated with $\vec{m}$ and $\vec{w}_2$ we have
\begin{align}\label{G}
    G_{\vec{m}} = \{\sum_{j=1}^q B_j w_{2,j} \mid
    B_j \mbox{ is a polynomial of degree } \leq m_j-1 \mbox{ for } j=1, \ldots, q \}
\end{align}
Hence $F_{\vec{n}}$ is the linear span of the functions in (\ref{spanF})
and $G_{\vec{m}}$ is the linear span of the functions in (\ref{spanG}).

Let $\vec{e}_k$ be the standard basis vector
\begin{align}
\vec{e}_k = (0,\ldots,0,1,0,\ldots,0), \quad \mbox{ where 1 is in
the $k$th position.}
\end{align}
We de not specify the length of the vector $\vec{e}_k$, but
this should be clear from the context.

\begin{lemma} \label{normality}
Suppose that $F_{\vec{n}}$ is $|\vec{n}|$-dimensional and let $|\vec{n}| = |\vec{m}|+1$.
\begin{enumerate}
\item[\rm (a)] Then $Q$ is a linear form {\rm (\ref{typeMform})} of
multiple orthogonal polynomials of mixed type for the pair
$(\vec{n}, \vec{m})$ if and only if
$ Q \in F_{\vec{n}} \cap G_{\vec{m}}^{\perp}$.
\item[\rm (b)] The pair $(\vec{n}, \vec{m})$ is a normal pair of indices
if and only if
$F_{\vec{n}} \cap G_{\vec{m}}^{\perp}$ is one-dimensional.
\item[\rm (c)] The pair $(\vec{n}, \vec{m})$ allows a type I normalization
with respect to the $k$th index of $\vec{m}$ if and only if
$F_{\vec{n}} \cap G_{\vec{m}+\vec{e}_k}^{\perp} = \{0\}$.
\item[\rm (d)] The pair $(\vec{n}, \vec{m})$ allows a type II normalization
with respect to the $k$th index of $\vec{n}$ if and only if
$F_{\vec{n}-\vec{e}_k} \cap G_{\vec{m}}^{\perp} = \{ 0 \}$.
\end{enumerate}
\end{lemma}

\begin{proof}
(a) This is immediate from the definitions.

(b) If $F_{\vec{n}} \cap G_{\vec{m}}^{\perp}$ is one-dimensional, then by part (a)
the linear form $Q$ is unique up to a multiplicative constant. Then the
polynomials $A_j$ are also unique up to a multiplicative constant, since
$F_{\vec{n}}$ is $|\vec{n}|$-dimensional. Thus $(\vec{n}, \vec{m})$ is
a normal pair of indices. The converse is obvious in view of part (a).

(c) Suppose $F_{\vec{n}} \cap G_{\vec{m} + \vec{e}_k}^{\perp} = \{0\}$.
Let $Q$ be a non-zero linear form of multiple orthogonal polynomials of
mixed type. Then $Q \in G_{\vec{m}}^{\perp}$ by part (a), but
$Q \not\in G_{\vec{m}+ \vec{e}_k}^{\perp}$.
Since $G_{\vec{m}} \subset G_{\vec{m} + \vec{e}_k}$ with codimension one,
it follows that $Q$ is not orthogonal to any function
in $G_{\vec{m} + \vec{e}_k} \setminus G_{\vec{m}}$ and since $x^{m_k} w_{2,k}$
belongs to this set, this implies
\[ \int Q(x) x^{m_k} w_{2,k}(x) dx \neq 0. \]
Then we can normalize $Q$ so that this integral is $1$ and
the pair $(\vec{n}, \vec{m})$ allows a type I normalization with
respect to the $k$th index.

On the other hand  if $F_{\vec{n}} \cap G_{\vec{m} + \vec{e}_k}^{\perp} \neq \{0\}$,
then any non-zero $Q$ in this space would be a non-zero linear form of
multiple orthogonal polynomials of mixed type such that
\[ \int Q(x) x^{m_k} w_{2,k}(x) dx  = 0. \]
Then either $(\vec{n}, \vec{m})$ is not normal, or if it is normal,
it does not allow a type I normalization with respect to the $k$th index.

(d)  Suppose $F_{\vec{n}-\vec{e}_k} \cap G_{\vec{m}}^{\perp} =
\{0\}$. Let $Q$ be a non-zero linear form of multiple orthogonal
polynomials of mixed type for the pair $(\vec{n},\vec{m})$. Then
$Q \in F_{\vec{n}}$ by part (a), but $Q \not\in
F_{\vec{n}-\vec{e}_k}$. This implies that $A_k$ has exact degree
$n_k-1$. Then we can normalize $Q$ so that the leading coefficient
of $A_k$ is $1$ and thus the pair $(\vec{n}, \vec{m})$ allows a
type II normalization with respect to the $k$th index.

On the other hand  if $F_{\vec{n}-\vec{e}_k} \cap
G_{\vec{m}}^{\perp} \neq \{0\}$, then any non-zero $Q$ in this
space would be a non-zero linear form of multiple orthogonal
polynomials of mixed type for the pair $(\vec{n},\vec{m})$ with
the degree of $A_k$ less than $n_k-1$. Then either $(\vec{n},
\vec{m})$ is not normal, or if it is normal, it does not allow a
type II normalization with respect to the $k$th index.
\end{proof}

We have the following easy corollary of Lemma \ref{normality}.

\begin{corollary}\label{normality-cor}
Suppose that
\begin{itemize}
\item[\rm (1)] $|\vec{n}| = |\vec{m}|$ and  $F_{\vec{n}}$ is a
$|\vec{n}|$-dimensional subspace of $L^2(\mathbb R)$. \item[\rm
(2)] $F_{\vec{n}} \cap G_{\vec{m}}^{\perp} = \{ 0 \}$.
\end{itemize}
Then the following hold:
\begin{enumerate}
\item[\rm (a)] For every $k=1, \ldots, q$ we have that $(\vec{n}, \vec{m}-\vec{e}_k)$
is a normal pair of multi-indices which allows a type I normalization with
respect to the $k$th index. Hence $Q^{(I,k)}_{\vec{n}, \vec{m}-\vec{e}_k}$
exists and is unique.
\item[\rm (b)] For every $k=1, \ldots, p$ we have that $(\vec{n} + \vec{e}_k, \vec{m})$
is a normal pair of multi-indices which allows a type II normalization with
respect to the $k$th index. Hence $Q^{(II,k)}_{\vec{n}+\vec{e}_k,\vec{m}}$
exists and is unique.
\end{enumerate}
\end{corollary}

\begin{proof}
Part (a) follows immediately from part (c) of Lemma \ref{normality} applied
to the multi-indices $\vec{n}$ and $\vec{m}-\vec{e}_k)$,
and part (b) follows from part (d) of Lemma \ref{normality} applied
to the multi-indices $\vec{n} + \vec{e}_k$ and $\vec{m}$.
\end{proof}

\section{The Riemann-Hilbert problem}
Fokas, Its, and Kitaev \cite{FIK} found a Riemann-Hilbert
problem that characterizes the orthogonal polynomials. Van Assche,
Geronimo, and Kuijlaars \cite{VAGK} extended this
Riemann-Hilbert problem to multiple orthogonal polynomials of type
I and type II. We are now going to give a further extension to
multiple orthogonal polynomials of mixed type.

\subsection{Riemann-Hilbert problem for multiple orthogonal polynomials
of mixed type}
Assume that $|\vec{n}| = |\vec{m}|$ and let $\vec{w}_1$ and
$\vec{w}_2$ be as before. For convenience we assume that $n_k> 0$ and
$m_l >0$ for $k=1,\ldots,p$ and $l=1,\ldots,q$. We can modify the
arguments in case one or several of the $n_k$ and $m_l$ are zero,
but we will not discuss that here.

Consider the
following Riemann-Hilbert problem: determine a $(p+q) \times
(p+q)$ matrix valued function $Y : \mathbb C \setminus \mathbb R \to
\mathbb C^{(p+q)\times(p+q)}$ such that
\begin{enumerate}
\item[(1)] $Y$ is analytic on $\mathbb{C} \setminus \mathbb{R}$,
\item[(2)]
for $ x \in \mathbb{R}$, we have
\begin{align}\label{RH1}
Y_{+}(x)=Y_{-}(x) \begin{bmatrix}
    I_p & W(x) \\
    0   & I_q
    \end{bmatrix} \end{align}
where $I_p$ and $I_q$ denote the identity matrices of sizes
$p$ and $q$, respectively, and
\begin{align}\label{defW}
W  = \vec{w}_1^{\;t}  \vec{w}_2
    = \begin{bmatrix}
    w_{1,1} w_{2,1} & w_{1,1} w_{2,2}  & \cdots &  w_{1,1} w_{2,q} \\
    w_{1,2} w_{2,1} & w_{1,2} w_{2,2}  & \cdots &  w_{1,2} w_{2,q} \\
    \vdots & \vdots & \ddots & \vdots  \\
    w_{1,p} w_{2,1} & w_{1,p}w_{2,2} & \cdots & w_{1,p} w_{2,q}
    \end{bmatrix}.
\end{align}
\item[(3)] as $z \to \infty$, we have that
\begin{align}\label{RH11}
Y(z)= \left(I_{p+q} + O\left(\frac{1}{z}\right) \right)
\begin{bmatrix}
z^{n_1} & 0 & 0 & \cdots & 0 & 0 & 0 \\
0 & z^{n_2} & 0 & \cdots & 0 & 0 & 0 \\
\vdots & \vdots & \ddots & \vdots & \vdots & \vdots & \vdots \\
0 & \cdots & 0 & z^{n_p} & 0 & \cdots & 0 \\
0 & \cdots & 0 & 0 & z^{-m_1} & \cdots & 0 \\
\vdots & \vdots & \vdots & \vdots & \vdots & \ddots & \vdots \\
0 & 0 & 0 & \cdots & 0 & 0 & z^{-m_q}
\end{bmatrix}
\end{align}
where $I_{p+q}$ denotes the identity matrix of size $p+q$.
\end{enumerate}

As before we use $F_{\vec{n}}$ and $G_{\vec{m}}$ to denote
the spaces (\ref{F}) and (\ref{G}) associated with the
indices $\vec{n}$ and $\vec{m}$.
The main result of this section is that the Riemann-Hilbert problem
has a unique solution if the conditions of Corollary \ref{normality-cor} are
satisfied.

\begin{theorem}\label{oplRH}
Let $|\vec{n}| = |\vec{m}|$. Suppose that $F_{\vec{n}}$ is
an $|\vec{n}|$-dimensional subspace of $L^2(\mathbb R)$ such
that $F_{\vec{n}} \cap G_{\vec{m}}^{\perp} = \{0\}$.
Then the above Riemann-Hilbert problem has a unique solution, given in terms of
the multiple orthogonal polynomials of mixed type with respect to
the vectors of weights $\vec{w}_1$ and $\vec{w}_2$. We have
\begin{itemize}
\item for $k,l = 1,\ldots,p$:
\begin{align} \label{block11}
    Y_{k,l}(z) =
A_{l,\vec{n}+\vec{e}_k,\vec{m}}^{(II,k)}(z),
\end{align}
 \item for $k=1,\ldots,q$ and $l=1,\ldots, p$:
\begin{align} \label{block21}
Y_{p+k,l}(z) = -2 \pi i
A_{l,\vec{n},\vec{m}-\vec{e}_k}^{(I,k)}(z),
\end{align}
\item for $k=1,\ldots,p$ and $l=1,\ldots,q$:
\begin{align} \label{block12}
Y_{k,p+l}(z) = \frac{1}{2 \pi i} \int
\frac{Q_{\vec{n}+\vec{e}_k,\vec{m}}^{(II,k)}(x) w_{2,l}(x)}{x-z} dx,
\end{align}
\item for  $k,l=1,\ldots,q$:
\begin{align} \label{block22}
Y_{p+k,p+l}(z) = - \int
\frac{Q_{\vec{n},\vec{m}-\vec{e}_k}^{(I,k)}(x)w_{2,l}(x)}{x-z} dx.
\end{align}
\end{itemize}
\end{theorem}
\begin{proof}
First note that the functions
$Q_{\vec{n}+\vec{e}_k,\vec{m}}^{(II,k)}$ and
$Q_{\vec{n},\vec{m}-\vec{e}_l}^{(I,l)}$ uniquely exist
for $k=1,\ldots,p$ and $l=1,\ldots,q$ by Corollary
\ref{normality-cor}. Partition $Y$ as
\begin{align}
Y(z) = \begin{bmatrix} K(z) & L(z) \\
    M(z) & N(z) \end{bmatrix},
\end{align}
where $K$ is a $p \times p$ matrix, $L$ is a $p \times q$ matrix, $M$ is
a $q \times p$ matrix, and $N$ is a $q \times q$ matrix.

The jump condition (\ref{RH1}) implies that
\begin{align}\label{analytic}
K_{k,l}^+(x) = K_{k,l}^-(x) \quad \mbox{ for } x \in \mathbb{R},
\end{align}
so that $K$ is analytic on the full complex plane.
From (\ref{analytic}), the asymptotic condition (\ref{RH11}), and
an extension of Liouville's theorem to polynomials, it then follows that
each diagonal element $K_{k,k}$, with $k=1,\ldots,p$, is
a monic polynomial of degree $n_k$, and that each off-diagonal
elements $K_{k,l}$ is a polynomial of degree at most $n_l-1$.

For the matrix $L$ the jump condition (\ref{RH1}) implies that for
$x \in \mathbb{R}$,
\begin{align}
L_{k,l}^+(x) = L_{k,l}^-(x) + Q_k(x) w_{2,l}(x) \quad \mbox{ for }
k=1,\ldots,p \mbox{ and } l = 1,\ldots,q,
\end{align}
where
\[ Q_k(x) = \sum_{j=1}^p K_{k,j}(x)w_{1,j}. \]
The Sokhotsky-Plemelj formula then gives that
\begin{align}\label{SP}
L_{k,l}(z)= \frac{1}{2\pi i}\int \frac{Q_k(x) w_{2,l}(x)}{x-z} dx
\quad  \mbox{ if } z \in \mathbb C \setminus \mathbb{R}.
\end{align}
If we now use the expansion
\begin{align}\label{exp}
\frac{1}{z-x} = \sum_{l=0}^{n-1} \frac{x^l}{z^{l+1}} +
\frac{x^n}{z^n}\frac{1}{z-x} \quad  \mbox{ for } n \in \mathbb{N},
\end{align}
we find that
\begin{align} \label{expandL}
L_{k,l}(z) = & - \sum_{j=0}^{n-1}  \frac{1}{ 2\pi i z^{j+1}}\int
    Q_k(x) x^j w_{2,l}(x) dx \nonumber
\\ & -\frac{1}{2\pi i z^n}\int \frac{Q_k(x) x^n
w_{2,l}(x)}{z-x} dx \quad \mbox{ for } k=1,\ldots,p \mbox{ and } l
= 1,\ldots,q.
\end{align}
The asymptotic condition (\ref{RH11}) gives that
\begin{align}
 \lim_{z \rightarrow \infty} L_{k,l}(z)z^{m_l} = 0,
\end{align}
such that from (\ref{expandL}) with $n = m_{l-1}$ we get
\begin{align}\label{Korth}
\int Q_k(x) x^j w_{2,l}(x) dx = 0 \quad \mbox{ for } j
=0,\ldots,m_l-1, k=1,\ldots,p \mbox{ and } l=1,\ldots,q.
\end{align}
Because $K_{k,l}$ is a polynomial of degree at most $n_l-1$ if $l
\neq k$ and $K_{k,k}$ is a monic polynomial of degree $n_k$, as
mentioned in the beginning of the proof, and because of
(\ref{Korth}) we see that
\begin{align}\label{K}
K_{k,l} = A_{l,\vec{n}+\vec{e}_k,\vec{m}}^{(II,k)}
\end{align}
for $k,l =1,\ldots,p$. Because of (\ref{SP}) and (\ref{K}) we see
that
\begin{align}
L_{k,l}(z) = \frac{1}{2 \pi i} \int
\frac{Q_{\vec{n}+\vec{e}_k,\vec{m}}^{(II,k)}(x)w_{2,l}(x)}{x-z}
dx \quad \mbox{ if } z \in \mathbb{C}\setminus  \mathbb{R}
\end{align}
for $k=1,\ldots,p$ and $l=1,\ldots,q$.
This proves the formulas (\ref{block11}) and (\ref{block12}).

\medskip
The jump condition (\ref{RH1}) shows that for $k=1, \ldots, q$ and $l = 1, \ldots, p$,
\begin{align}
M_{k,l}^+(x) = M_{k,l}^-(x) \quad \mbox{ for } x \in \mathbb{R},
\end{align}
and consequently $M$ is also analytic on the full complex plane.
In the same way as for the matrix $K$, the asymptotic condition
(\ref{RH11}) implies that each $M_{k,l}$ is a polynomial of degree
$\leq n_l - 1$. The jump condition (\ref{RH1}) also implies that

\begin{align}
N_{k,l}^+(x) = N_{k,l}^-(x) + Q_k(x) w_{2,l}(x) \quad \mbox{ for }
k,l=1,\ldots,q \mbox{ and } x \in \mathbb{R},
\end{align}
with
\[ Q_k(x) = \sum_{j=1}^p M_{k,j}(x) w_{1,j}(x). \]
Using the Sokhotsky-Plemelj formula we get that
\begin{align}\label{condN}
N_{k,l}(z) = \frac{1}{2\pi i} \int \frac{Q_k(x) w_{2,l}(x)}{x-z}dx
\quad \mbox{ for } k,l=1,\ldots,q \mbox{ if } z \in \mathbb{C}
\setminus \mathbb{R}.
\end{align}
The asymptotic condition (\ref{RH11}) gives that
\[\begin{cases}
\lim\limits_{z \to \infty} N_{k,l}(z) z^{m_l} = 0 \quad \mbox{ if } k \neq l, \\[10pt]
\lim\limits_{z \to \infty} N_{k,k}(z) z^{m_k} = 1, \\
\end{cases}\]
Using these conditions and the expansion (\ref{exp})  we get that
\begin{align}\label{testje}
\begin{cases} \int Q_k(x) x^j w_{2,l}(x)dx = 0 \quad \mbox{
for } j=0,1,\ldots, m_l - 1 \mbox{ if } k \neq l, \\
 \int Q_k(x) x^j w_{2,k}(x)dx = 0 \quad
\mbox{ for } j=0,1,\ldots, m_k - 2,
\\
 \int Q_k(x) x^{m_k-1} w_{2,k}(x)dx = -2\pi i.
\end{cases}
\end{align}
The degree of the polynomials $M_{k,l}$ and the orthogonality
conditions (\ref{testje}) imply that $M_{k,l} = -2 \pi i
A_{l,\vec{n},\vec{m}-\vec{e}_k}^{(I,k)}$ for $k =
1,\ldots,q$ and $l=1,\ldots,p$. Because of (\ref{condN}) we get
that
\begin{align}
N_{k,l}(z) = -\int \frac{Q_{\vec{n},\vec{m}-\vec{e}_k}^{(I,k)}(x) w_{2,l}(x)}{x-z}dx
\end{align}
for $k,l=1,\ldots,q$ and $z \in \mathbb{C} \backslash \mathbb{R}$.
This proves the formulas (\ref{block21}) and (\ref{block22}). This
completes the proof of Theorem \ref{oplRH}.
\end{proof}

\subsection{Riemann-Hilbert problem for the inverse}
By standard arguments it follows that $\det Y(z) \equiv 1$
for $z \in \mathbb C \setminus \mathbb R$ so that the
inverse $Y^{-1}(z)$ exists and is analytic for
$z \in \mathbb C \setminus \mathbb R$.  Define
\begin{align}\label{Y^(-t)}
    X(z)=Y^{-t}(z).
\end{align}
From the Riemann-Hilbert problem for $Y$ it is then straightforward
to check that $X$ is the solution of the following $(p+q) \times (p+q)$
matrix valued Riemann-Hilbert problem:
\begin{enumerate}
\item[\rm (1)] $X : \mathbb C \setminus \mathbb R \to
    \mathbb C^{(p+q)\times (p+q)}$ is analytic,
\item[\rm (2)]
for $ x \in \mathbb{R}$, we have
\begin{align}\label{jump}
  X_{+}(x)=X_{-}(x) \begin{bmatrix} I_p & 0 \\
    -W^{t}(x)   & I_q \end{bmatrix}
\end{align}
where  $W$ is given by (\ref{defW}),
\item[\rm (3)] as $z \to \infty$, we have that
\begin{align}\label{asym}
X(z)= \left(I_{p+q} + O\left(\frac{1}{z}\right) \right)
    \begin{bmatrix}
    z^{-n_1} & 0 & 0 & \cdots & 0 & 0 & 0 \\
    0 & z^{-n_2} & 0 & \cdots & 0 & 0 & 0 \\
    \vdots & \vdots & \ddots & \vdots & \vdots & \vdots & \vdots \\
    0 & \cdots & 0 & z^{-n_p} & 0 & \cdots & 0 \\
    0 & \cdots & 0 & 0 & z^{m_1} & \cdots & 0 \\
    \vdots & \vdots & \vdots & \vdots & \vdots & \ddots & \vdots \\
    0 & 0 & 0 & \cdots & 0 & 0 & z^{m_q}
    \end{bmatrix}.
\end{align}
\end{enumerate}

The solution of the Riemann-Hilbert problem for $X$
can again be written
in terms of multiple orthogonal polynomials of mixed type, but
with the roles of the vectors of weights $\vec{w}_1$ and
$\vec{w}_2$ as well as the multi-indices $\vec{n}$ and $\vec{m}$
interchanged. Therefore we use the full notation
$Q_{\vec{m},\vec{n}}(x;\vec{w}_2,\vec{w}_1)$.

\begin{lemma}\label{oplRH2}
Suppose that the conditions of Theorem {\rm \ref{oplRH}} are
satisfied. Then the solution of the above Riemann-Hilbert problem
 has a unique solution given by
\begin{itemize}
\item For $k,l=1,\ldots,p$:
\begin{align} \label{solRH11}
X_{k,l}(z) = -\int
\frac{Q_{\vec{m},\vec{n}-\vec{e}_k}^{(I,k)}(x;\vec{w}_2,\vec{w}_1)w_{1,l}(x)}{x-z}dx,
\end{align}
\item for $k=1,\ldots,q$ and $l=1,\ldots,p$:
\begin{align} \label{solRH12}
X_{p+k,l}(z) = - \frac{1}{2\pi i}\int
\frac{Q_{\vec{m}+\vec{e}_k,\vec{n}}^{(II,k)}(x;\vec{w}_2,\vec{w}_1)w_{1,l}(x)}{x-z}dx,
\end{align}
\item for $k=1,\ldots,p$ and $l=1,\ldots,q$:
\begin{align}\label{solRH21}
 X_{k,p+l}(z) = 2 \pi i A_{l,\vec{m},\vec{n}-\vec{e}_k}^{(I,k)}(z;\vec{w}_2,\vec{w}_1),
\end{align}
\item for $k,l=1,\ldots,q$:
\begin{align}\label{solRH22}
X_{p+k,p+l}(z) =
A_{l,\vec{m}+\vec{e}_k,\vec{n}}^{(II,k)}(z;\vec{w}_2,\vec{w}_1).
\end{align}
\end{itemize}
\end{lemma}
\begin{proof}
The lemma can be proven in the same way as Theorem \ref{oplRH},
but it is also possible to derive it directly from Theorem
\ref{oplRH} as follows. Let $U$ be the solution of the Riemann-Hilbert
problem described in Theorem \ref{oplRH}, but with the roles of
the vectors of weights $\vec{w}_1$ and $\vec{w}_2$ as well as the
multi-indices $\vec{n}$ and $\vec{m}$ interchanged. By comparing
the jump conditions and the asymptotic conditions of the
Riemann-Hilbert problem of $X$ and $U$, we can easily see that
\begin{align}
U(z) = \begin{bmatrix} 0 & -I_q \\ I_p & 0 \end{bmatrix}
X(z)
\begin{bmatrix} 0 & I_p\\ -I_q & 0 \end{bmatrix}.
\end{align}
Theorem \ref{oplRH} therefore implies that the solution of the
Riemann-Hilbert problem defined above is unique and is given by
the formulas (\ref{solRH11})--(\ref{solRH22}).
\end{proof}

\section{The kernel}
Suppose that $\vec{w}_1 = (w_{1,1}, w_{1,2}, \ldots, w_{1,p})$ and
$\vec{w}_2 = (w_{2,1}, w_{2,2}, \ldots, w_{2,q})$ are two vectors
of weights on the real line, and define $F_{\vec{n}}$ and
$G_{\vec{m}}$ as in (\ref{F}) and (\ref{G}). Suppose as in section
3 that $|\vec{n}| =|\vec{m}|=n$ and that $F_{\vec{n}}$ and
$G_{\vec{m}}$ are both $n$-dimensional subspaces of
$L^2(\mathbb{R})$. Two bases $\phi_1, \ldots, \phi_n$ of
$F_{\vec{n}}$ and  $\psi_1, \ldots, \psi_n$ of $G_{\vec{m}}$ are
called biorthogonal if
\[ \int \phi_j(x) \psi_k(x) dx = \delta_{j,k}. \]
The following lemmas are well-known, but we include their proofs
for completeness.
\begin{lemma} \label{bases}
There exist biorthogonal bases for $F_{\vec{n}}$ and $G_{\vec{m}}$
if and only if $F_{\vec{n}} \cap G_{\vec{m}}^{\bot} = \{0\}$.
\end{lemma}
\begin{proof}
First suppose that $\phi_1,\ldots,\phi_n \in F_{\vec{n}}$ and
$\psi_1,\ldots,\psi_n \in G_{\vec{m}}$ are biorthogonal bases.
Then every $f \in F_{\vec{n}}$ can be written as
\begin{align}
    f  = \sum_{j=1}^n c_j \phi_j
\end{align}
with
\begin{align}
    c_j = \int f(x) \psi_j(x) dx, \qquad j = 1, \ldots, m.
\end{align}
Since every $c_j = 0$ if $f \in G_{\vec{m}}^{\bot}$ it follows
that $F_{\vec{n}} \cap G_{\vec{m}}^{\bot} = \{0\}$.

Conversely, suppose that $F_{\vec{n}} \cap G_{\vec{m}}^{\bot} = \{0\}$.
Let $\psi_1,\ldots,\psi_n$ be any basis of $G_{\vec{n}}$. Consider the linear mapping
\begin{align}
\mathcal{F}: F_{\vec{n}} \to \mathbb{R}^n: f \mapsto (\int
f(x)\psi_k(x) dx)_{k=1,\ldots,n}.
\end{align}
Because $F_{\vec{n}} \cap G_{\vec{m}}^{\bot} = \{0\}$ it is clear
that $\mathcal{F}$ is injective. Since the
dimensions of $\mathbb{R}^n$ and $F_{\vec{n}}$ are equal, and the mapping is linear,
$\mathcal{F}$ is bijective. Consequently there exist functions
$\phi_j \in F_{\vec{n}}$, $j=1,\ldots,n$ such that $\int \phi_j(x)
\psi_k(x) dx = \delta_{j,k}$. Then $\phi_1, \ldots, \phi_n$ is a basis
of $F_{\vec{n}}$ which is biorthogonal to $\psi_n, \ldots, \psi_n$.
\end{proof}

From now on, we  assume that
$|\vec{n}| = |\vec{m}|$ and $F_{\vec{n}} \cap G_{\vec{m}}^{\bot} = \{0\}$.
According to Lemma \ref{bases}
there exist biorthogonal bases $\phi_1,\ldots,\phi_{n}$ of
$F_{\vec{n}}$ and $ \psi_1,\ldots,\psi_n$ of $G_{\vec{m}}$.
We define the kernel $K(x,y)$ as
\begin{align}\label{defK}
K(x,y) = \sum_{j=1}^n \phi_j(x)\psi_j(y).
\end{align}

\begin{lemma}\label{operator}
The kernel $K(x,y)$ is the kernel of the (non-orthogonal) projection operator
onto $F_{\vec{n}}$ parallel to $G_{\vec{m}}^{\bot}$.
\end{lemma}
\begin{proof}
Define the operator $K$ on $L^2(\mathbb{R})$ as
\begin{align}
(K h)(x) = \int K (x,y)h(y)dy.
\end{align}
By the definition (\ref{defK}) of the function $K(x,y)$ it is
clear that $K h = 0$ if $h \in G_{\vec{m}}^{\bot}$. In the same
way it is obvious that $K h=h$ if $h=\phi_k$, for $k=1,\ldots,n$.
By linearity it then follows that $K h=h$ for every $h \in
F_{\vec{n}}$. Because $F_{\vec{n}} \oplus G_{\vec{m}}^{\bot}=
L^2(\mathbb{R})$, Lemma \ref{operator} follows immediately.
\end{proof}
Lemma \ref{operator} also implies that the kernel $K$ is
independent of the chosen biorthogonal bases. Now we arrive at the
main result of this paper:

\begin{theorem}\label{theorem}
The kernel $K (x,y)$, as defined in {\rm(\ref{defK})}, can be
written in terms of the solution of the Riemann-Hilbert problem
for $Y$ of section 3.1 in the following way:
\begin{align} \label{defKinY}
K (x,y)
 & = \frac{1}{2 \pi i(x-y)}
  \begin{bmatrix}
0 & \cdots & 0 & w_{2,1}(y) & \cdots & w_{2,q}(y)
 \end{bmatrix}
 Y_+^{-1}(y) Y_+(x)
 \begin{bmatrix} w_{1,1}(x) \\ \vdots \\ w_{1,p}(x)\\ 0\\
\vdots \\ 0 \end{bmatrix}.
\end{align}
\end{theorem}

\begin{proof}
Define the operator $L$ on $L^2(\mathbb{R})$ as
\begin{align}
(Lh)(x) = \int L(x,y)h(y)dy,
\end{align}
where $L(x,y)$ denotes the right-hand side of (\ref{defKinY}).
It is enough to prove the following two things:
\begin{enumerate}
\item[(a)] $Lh = 0$ if $h \in G_{\vec{m}}^{\bot}$,
\item[(b)] $Lh = h$ if $h \in F_{\vec{n}}$.
\end{enumerate}
Indeed, if this is the case, then $L$ is the projection operator
onto $F_{\vec{n}}$ parallel to $G_{\vec{m}}^{\bot}$, and according to
Lemma \ref{operator}, we get that $K(x,y)$ is the kernel of $L$,
and consequently $K(x,y) = L(x,y)$.

\medskip

(a) First let $h \in G_{\vec{m}}^{\bot}$. For ease of notation
we will use $\begin{bmatrix}  \vec{w}_1(x) & \vec{0}
\end{bmatrix}$ and $\begin{bmatrix} \vec{0} & \vec{w}_2(y) \end{bmatrix}$ instead of
$\begin{bmatrix} w_{1,1}(x) & \cdots & w_{1,p}(x) & 0 & \cdots &  0
\end{bmatrix} $ and $\begin{bmatrix} 0&
\cdots &  0& w_{2,1}(x) & \cdots & w_{2,q}(x)
\end{bmatrix}$, where in each case the number of zeros is such that
the length of the vectors is $p+q$. We then have that
\begin{align}\label{opsplsom}
(Lh)(x)  =  & \frac{1}{2\pi i}\int h(y) \begin{bmatrix} \vec{0} & \vec{w}_2(y)
\end{bmatrix} \frac{Y_+^{-1}(y)-Y_+^{-1}(x)}{x-y}Y_+(x)
\begin{bmatrix}  \vec{w}_{1}(x) & \vec{0}  \end{bmatrix}^t dy \nonumber \\
&  + \frac{1}{2 \pi i}\int h(y) \begin{bmatrix} 0 & \vec{w}_2(y)
\end{bmatrix} \frac{Y_+^{-1}(x)}{x-y}Y_+(x)
\begin{bmatrix}  \vec{w}_{1}(x) & \vec{0}  \end{bmatrix}^t dy.
\end{align}
Because $Y_+^{-1}(x)Y_+(x) = I$ and $\begin{bmatrix} \vec{0} &
\vec{w}_2(y)\end{bmatrix} \begin{bmatrix}  \vec{w}_{1}(x) &
\vec{0} \end{bmatrix}^t  = 0$, the second term of (\ref{opsplsom})
is equal to zero. The form of the solution of the  Riemann-Hilbert
problem for $X = Y^{-t}$ as given by (\ref{solRH21}) and (\ref{solRH22})
implies that the last $q$ rows of
\begin{align}
\frac{Y_+^{-1}(y)-Y_+^{-1}(x)}{x-y}
\end{align}
consist of polynomials in the variable $y$ such that
for $j =1, \ldots, q$, and $k=1,\ldots, p+q$,
\begin{align}
\deg \left[\frac{Y_+^{-1}(y)-Y_+^{-1}(x)}{x-y}\right]_{p+j,k} \leq m_j - 1.
\end{align}
This implies that for each fixed $x \in \mathbb R$,
each entry of the row vector
\[
\begin{bmatrix} \vec{0} & \vec{w}_2(y)
\end{bmatrix} \frac{Y_+^{-1}(y)-Y_+^{-1}(x)}{x-y}
\]
belongs to $G_{\vec{m}}$. Because $h \in G_{\vec{m}}^{\bot}$,
the first term of (\ref{opsplsom}) is equal to zero
as well. Thus $Lh = 0$ and this proves (a).

\medskip

(b) Now let $h \in F_{\vec{n}}$. Then $h=\sum_{j=1}^p A_j(x)w_{1,j}(x)$, where
$A_j$ is a polynomial of degree less than or equal to $n_j-1$. We write
$\vec{A} = (A_1, \ldots, A_p)$ and
$\begin{bmatrix} \vec{A}(x) & \vec{0} \end{bmatrix} =
\begin{bmatrix} A_1(x) & \cdots & A_p(x) & 0 & \cdots
& 0 \end{bmatrix}$. We then have that $h(x) = \begin{bmatrix} \vec{A}(x) & \vec{0} \end{bmatrix}
    \begin{bmatrix} \vec{w}_1(x) & 0 \end{bmatrix}^t$ and so
\begin{align}
(Lh)(x) = & \frac{1}{2 \pi i} \int
\frac{\begin{bmatrix} \vec{A}(y)-\vec{A}(x) &  \vec{0} \end{bmatrix} }{x-y}
\begin{bmatrix} \vec{w}_1(y) & \vec{0} \end{bmatrix}^t
\begin{bmatrix} \vec{0} & \vec{w}_2(y)\end{bmatrix} Y_+^{-1}(y) Y_+(x)
\begin{bmatrix} \vec{w}_{1}(x) & \vec{0} \end{bmatrix}^t dy \nonumber \\[10pt]
& + \frac{1}{2\pi i} \int \begin{bmatrix} \vec{A}(x) & \vec{0} \end{bmatrix}
 \begin{bmatrix} \vec{w}_{1}(y) & \vec{0} \end{bmatrix}^t
\begin{bmatrix} \vec{0} & \vec{w}_2(y)\end{bmatrix} Y_+^{-1}(y) \frac{Y_+(x)}{x-y}
\begin{bmatrix} \vec{w}_{1}(x) & \vec{0} \end{bmatrix}^t dy. \label{int}
\end{align}

We will deal first with the first term in the right-hand side of (\ref{int}).
We have the combination $\begin{bmatrix} \vec{0} & \vec{w}_2(y) \end{bmatrix} Y_+^{-1}(y)$
which is a row vector whose $k$th entry is
\[ \sum_{l=1}^q \left(Y_+^{-1}\right)_{p+l,k}(y) w_{2,l}(y) =
    \sum_{l=1}^q X_{k,p+l}(y) w_{2,l}(y) \]
since $Y^{-1} = X^t$, see (\ref{Y^(-t)}). The functions
$X_{k,p+l}$ are certain multiple orthogonal polynomials of mixed
type given explicitly by formulas (\ref{solRH21}) and
(\ref{solRH22}). Then it follows that
\begin{align} \label{kthentry1}
    \left( \begin{bmatrix} \vec{0} & \vec{w}_2(y) \end{bmatrix} Y_+^{-1}(y) \right)_k
    = 2\pi i Q_{\vec{m},\vec{n}-\vec{e}_k}^{(I,k)}(y;\vec{w}_2,\vec{w}_1)
    \qquad \mbox{for }k = 1,\ldots,p
     \end{align}
and
\begin{align} \label{kthentry2}
   \left( \begin{bmatrix} \vec{0} & \vec{w}_2(y) \end{bmatrix} Y_+^{-1}(y) \right)_{p+k}
     = Q_{\vec{m}+\vec{e}_{k}, \vec{n}}^{(II,k)}(y; \vec{w}_2, \vec{w}_1)
   \qquad \mbox{for } k =  1, \ldots, q.
    \end{align}
Since $\frac{A_j(y)-A_j(x)}{x-y}$ is a polynomial of degree
$\leq n_j-2$ in the variable $y$ for $j=1,\ldots,p$,
we have for each fixed $x \in \mathbb R$,
\begin{align}
\frac{\begin{bmatrix} \vec{A}(y)-\vec{A}(x) &  \vec{0}
\end{bmatrix}}{x-y} \begin{bmatrix} \vec{w}_1(y) & \vec{0}
\end{bmatrix}^t \in F_{\vec{n}-\sum_{j=1}^p \vec{e}_j}.
\end{align}
From the defining properties of the multiple orthogonal polynomials
of mixed type, it follows that each of the functions (\ref{kthentry1})
and (\ref{kthentry2}) is orthogonal to $F_{\vec{n}-\sum_{j=1}^p \vec{e}_j}$.
Then it follows that the first integral in the right-hand side
of (\ref{int}) is zero for every $x$.
\medskip

Now we come to the second term in the right-hand side of (\ref{int}).
We are going to show that for every $x \in \mathbb R$,
\begin{align} \label{toprove}
    \frac{1}{2\pi i} \int \begin{bmatrix} \vec{w}_1(y) & \vec{0} \end{bmatrix}^t
\begin{bmatrix} \vec{0} & \vec{w}_2(y) \end{bmatrix} Y_+^{-1}(y) \frac{Y_+(x)}{x-y} dy
= \begin{bmatrix} I_p & * \\
    * & * \end{bmatrix},
\end{align}
where $*$ represents an unspecified unimportant entry (which may actually be
a divergent integral).  Having (\ref{toprove}) we easily see that the second term in
the right-hand side of (\ref{int}) reduces to $\begin{bmatrix} \vec{A}(x) & \vec{0} \end{bmatrix}
\begin{bmatrix} \vec{w}_1(x) & \vec{0} \end{bmatrix}^t = h(x)$,
independent of what the unspecified entries are (even if they are
divergent integrals).

In order to establish (\ref{toprove}) we note that the
jump condition (\ref{RH1}) written in the form
\[
Y_+(y) = Y_-(y)\left(I + \begin{bmatrix}  \vec{w}_{1}(y) & \vec{0}
\end{bmatrix}^t \begin{bmatrix}  \vec{0} & \vec{w}_2(y)
\end{bmatrix}\right)
\]
implies that
\begin{align}
\begin{bmatrix} \vec{w}_1(y) & \vec{0} \end{bmatrix}^t
\begin{bmatrix}  \vec{0} & \vec{w}_2(y)\end{bmatrix}
Y_+^{-1}(y) & = Y_-^{-1}(y)\left(Y_+(y) -
Y_-(y)\right) Y_+^{-1}(y) \nonumber \\[10pt]  &
= Y_-^{-1}(y) -
Y_+^{-1}(y). \label{verschil}
\end{align}
Thus the left-hand side of (\ref{toprove}) is
\begin{align}
\frac{1}{2 \pi i} \int \frac{Y_-^{-1}(y) - Y_+^{-1}(y)}{x-y} Y_+(x) dy
\end{align}
Let $z \in \mathbb C$ with $\Im z > 0$. Then
$Y^{-1}(y)/(z-y)$ is analytic in the lower half plane and
from the Riemann-Hilbert problem satisfied by $X = Y^{-t}$,
it follows that for $k,=1, \ldots, p$, $l=1, \ldots, p+q$,
\begin{align} \label{decay}
    \frac{\left[Y^{-1}(y)\right]_{k,l}}{z-y}  = O(y^{-n_k-1}) \qquad
    \mbox{as } y \to \infty.
\end{align}
It follows that (here we use $n_k \geq 1$)
\begin{align} \label{partsum1}
    \frac{1}{2\pi i} \int \frac{\left[Y_-^{-1}(y)\right]_{k,l}}{z-y} dy = 0,
        \qquad \mbox{for } k =1, \ldots, p, \quad l=1, \ldots, p+q.
\end{align}
Similarly, we have that $Y^{-1}(y)/(z-y)$ is analytic in the
upper half plane but with a pole at $y=z$. Then if we calculate
the same integral as in (\ref{partsum1}) but with $Y_-$ replaced
by $Y_+$, and we use the decay property (\ref{decay}),
the only contribution comes form the residue at $y=z$ and the result is
\begin{align} \label{partsum2}
    \frac{1}{2\pi i} \int \frac{\left[Y_+^{-1}(y) \right]_{k,l}}{z-y} dy =
    - \left( Y^{-1}(z) \right)_{k,l},
    \qquad \mbox{for } k =1, \ldots, p, \quad l=1, \ldots, p+q.
\end{align}
From (\ref{partsum1}) and (\ref{partsum2}) it follows that
\begin{align}
\left(\frac{1}{2 \pi i} \int \frac{Y_-^{-1}(y) - Y_+^{-1}(y)}{z-y}dy \right) Y(z)
    = \begin{bmatrix} I_p & * \\ * & * \end{bmatrix}
        \qquad \mbox{for } \Im z > 0.
\end{align}
Letting $z \to x \in \mathbb R$, it follows that
\begin{align}
\frac{1}{2 \pi i} \int \frac{Y_-^{-1}(y) - Y_+^{-1}(y)}{x-y} Y_+(x) dy
    = \begin{bmatrix} I_p & * \\ * & * \end{bmatrix}
        \qquad \mbox{for } x \in \mathbb R,
\end{align}
which implies (\ref{toprove}) by (\ref{verschil}).
As noted after (\ref{toprove}) it then follows that the
second term in the right-hand side of (\ref{int})
is equal to $h(x)$. Since we already know that the
first term is equal to $0$, we have proved that
$Lh = h$. This completes the proof of Theorem \ref{theorem}.
\end{proof}

\section{The Christoffel-Darboux formula}

Theorem \ref{theorem} implies a  Christoffel-Darboux
formula for multiple orthogonal polynomials of mixed type.
We assume as before that $|\vec{n}| = |\vec{m}|$ and that
$F_{\vec{n}}$ and $G_{\vec{m}}$ are both $|\vec{n}|$-dimensional
subspaces of $L^2(\mathbb R)$ such that $F_{\vec{n}} \cap G_{\vec{m}}^{\bot}
= \{ 0 \}$.
\begin{corollary}\label{christoffelke}
Let $K $ be the  kernel defined in {\rm (\ref{defK})}. We can
write the kernel in terms of the multiple orthogonal
polynomials of mixed type defined in {\rm (\ref{normI})} and
{\rm(\ref{normII})} as follows:
\begin{align}\label{christoffel}
(x-y)K (x,y)  = &  \sum_{j=1}^{p}
Q_{\vec{n}+\vec{e}_j,\vec{m}}^{(II,j)}(x;\vec{w}_1,\vec{w}_2)
Q_{\vec{m},\vec{n}-\vec{e}_j}^{(I,j)}(y;\vec{w}_2,\vec{w}_1) \nonumber \\ & -
\sum_{k=1}^q
Q_{\vec{n},\vec{m}-\vec{e}_k}^{(I,k)}(x;\vec{w}_1,\vec{w}_2)
Q_{\vec{m}+\vec{e}_k,\vec{n}}^{(II,k)}(y;\vec{w}_2,\vec{w}_1).
\end{align}
\end{corollary}
\begin{proof}
The kernel $K $ and the multiple orthogonal polynomials of mixed
type are well defined because of Corollary \ref{normality-cor} and
Lemma \ref{bases}. The entries of $Y_+(x)$ in the first $p$
columns are given by (\ref{block11}) and (\ref{block21}). The
entries of $Y_+^{-1}(y) = X^t(y)$ in the last $q$ rows are given
by (\ref{solRH21}) and (\ref{solRH22}). Inserting these formulas
into (\ref{defKinY}) we arrive at (\ref{christoffel}).
\end{proof}

\begin{remark}
The usual monic orthogonal polynomials on the real line with
weight function $w(x)$ satisfy the classical Christoffel-Darboux formula
(\ref{christoffel0}).
By putting $p=q=1$ in formula (\ref{christoffel}) and taking into
account that the type II normalization for the multiple
orthogonal polynomials of mixed type is different from the
normalization used for monic orthogonal polynomials, we can
see that (\ref{christoffel}) reduces to (\ref{christoffel0})
in case $p=q=1$.

In \cite{BK1} the special case $p=1$ and $q=2$ is considered
in connection with random matrices with external source.
This leads to a kernel built out of
multiple orthogonal polynomials with respect to 2 different
weights for which a Christoffel-Darboux kernel was given.
In \cite{DK} the Christoffel-Darboux formula was generalized
to multiple orthogonal polynomials with respect to $q$ different weights:
\begin{align}\label{cd}
    (x-y)K (x,y)  & =
    P_{\vec{n}}(x)Q_{\vec{n}}(y)-\sum_{k=1}^{q}
    \frac{h_{\vec{n}}^{(k)}}{h_{\vec{n} -\vec{e}_{k}}^{(k)}}
    P_{\vec{n}-\vec{e}_k}(x) Q_{\vec{n}+\vec{e}_k}(y).
\end{align}
Here $P_{\vec{n}}$ is the multiple orthogonal polynomial of type
II, $Q_{\vec{n}}$ is the linear form constructed out of the
multiple orthogonal polynomials of type I, with the type I
normalization as described in section 2, and the
$h_{\vec{n}}^{(k)}$ and $h_{\vec{n}-\vec{e}_k}^{(k)}$
are certain constants. The formula (\ref{cd}) is the
special case $p=1$ of (\ref{christoffel}).
\end{remark}
\begin{remark}
When we take the multi-indices $\vec{n}$ and $\vec{m}$ in a way such
that
\begin{align}
n_1 \leq n_2 \leq \cdots \leq n_p \leq n_1+1
\quad \mbox{ and } \quad m_1 \leq
m_2 \leq \cdots \leq m_q \leq m_1+1,
\end{align}
 then the multiple orthogonal polynomials of mixed type are
vector polynomials orthogonal with respect to the weight
\begin{align}
W(x) & =
\begin{bmatrix}
w_{1,1}(x)w_{2,1}(x) & w_{1,1}(x)w_{2,2}(x) & \cdots &
w_{1,1}(x)w_{2,q}(x) \\
w_{1,2}(x)w_{2,1}(x) & w_{1,2}(x)w_{2,2}(x) & \cdots &
w_{1,2}(x)w_{2,q}(x) \\
\vdots & \vdots & \vdots & \vdots  \\
w_{1,p}(x)w_{2,1}(x) & w_{1,p}(x)w_{2,2}(x) & \cdots &
w_{1,p}(x)w_{2,q}(x)
\end{bmatrix}.
\end{align}
In \cite{SVI}, Sorokin and Van Iseghem obtained a
Christoffel-Darboux formula for vector polynomials. Their
formula has $\frac{p(p+1)}{2} + \frac{q(q+1)}{2}$ terms,
while ours has only $p+q$ terms.
\end{remark}

\section{Non intersecting Brownian motions}

Our motivation for introducing the multiple orthogonal
polynomials of mixed type came from the theory of
non-intersecting Brownian motions.
Consider $n$ one-dimensional Brownian motions which start at $n$
fixed points $a_1 < a_2 < \cdots < a_n$ at time $t=0$ and end at
$n$ fixed points $b_1 < b_2< \cdots < b_n$ at time $t=1$. Let
$p_{n,t}(x_1,\ldots,x_n)$ denote the probability density that at
time $t$, with $0 < t < 1$, the paths are at the positions $x_1,
\ldots, x_n$, conditioned on the event that the paths do not
intersect in the full time interval $(0,1)$. Then it follows from
a result of Karlin and McGregor \cite{KMcG} that
\begin{align}\label{karlin}
    p_{n,t}(x_1,\ldots,x_n) = \frac{1}{Z_n}
    \det(P(t,a_j,x_k))_{j,k=1}^n \det(P(1-t,b_j,x_k))_{j,k=1}^n,
\end{align}
where
\begin{align} \label{transprob}
P(t,a,x) = \frac{1}{\sqrt{2\pi t}} e^{-\frac{1}{2t}(x-a)^2}
\end{align}
is the transition probability for the one-dimensional Brownian
motion and $Z_n$ is a normalization constant. Note that
(\ref{karlin}) is an example of a biorthogonal ensemble \cite{Bor}.

Consider now the confluent case that some of the starting points
and some of the end points coincide. Suppose that the $n$
non-intersecting Brownian motions start at $p$ different points
$a_j$, $j=1,\ldots,p$, where $a_j$ appears with multiplicity $n_j$, and
end at $q$ different points $b_j$, $j=1,\ldots,q$, where $b_j$ appears
with multiplicity $m_j$. Let $\vec{n} = (n_1, \ldots,n_p)$
and $\vec{m} = (m_1, \ldots, m_q)$ and
\begin{align}
F_{\vec{n}} = \{\sum_{j=1}^p A_j(x) P(t,a_j,x) \mid A_j \mbox{
polynomial with } \deg(A_j) \leq n_j-1 \}
\end{align}
and
\begin{align}
G_{\vec{m}} = \{\sum_{j=1}^q B_j(x)P(1-t,b_j,x) \mid B_j \mbox{
polynomial with } \deg(B_j) \leq m_j-1  \}.
\end{align}
So these are the spaces (\ref{F}) and (\ref{G}) associated with the
vectors of weights $\vec{w}_1 = (w_{1,1}, \ldots, w_{1,p})$
and $\vec{w}_2 = (w_{2,1}, \ldots, w_{2,p})$ where
\[ w_{1,j}(x) = P(t,a_j,x), \qquad \mbox{ for } j=1, \ldots, p, \]
and
\[ w_{2,j}(x) = P(1-t,b_j,x), \qquad \mbox{ for } j=1, \ldots, q. \]

\begin{lemma} \label{Cheblemma}
The spaces $F_{\vec{n}}$ and $G_{\vec{m}}$ are $n$-dimensional and
\begin{align} \label{Chebsys}
    F_{\vec{n}} \cap G_{\vec{m}}^{\bot} = \{0\}.
\end{align}
\end{lemma}
\begin{proof}
The statement about the dimensions are obvious. The proof of
(\ref{Chebsys}) is based on the following facts:
\begin{enumerate}
\item[(1)] Both $F_{\vec{n}}$ and $G_{\vec{m}}$ are Chebyshev
spaces on $\mathbb R$, which means that any non-zero function in
one of these spaces has at most $n-1$ zeros on $\mathbb R$, see e.g. \cite{BE}.
To show this, we note that by an example given in  \cite[Chapter 4, \S4]{NS},
the functions
\begin{align} \label{chebbasis}
    e^{\beta_1x},\ldots,x^{n_1-1}e^{\beta_1x},\ldots,e^{\beta_px},\ldots,x^{n_p-1}e^{\beta_px}
\end{align}
form a Chebyshev system of order $n-1$ on $\mathbb{R}$ whenever $\beta_1,\ldots,\beta_p$ are
distinct real numbers. Taking $\beta_j = \frac{1}{2t} a_j$, and multiplying
the functions (\ref{chebbasis}) by the common factor $e^{-\frac{1}{2t} x^2}$, we obtain
a basis of $F_{\vec{n}}$, and so $F_{\vec{n}}$ is a Chebyshev space on $\mathbb R$. Similarly
we have that $G_{\vec{m}}$ is a Chebyshev spaces on $\mathbb R$.
\item[(2)] For
any set of distinct real points $x_1, \ldots, x_m$ with $m \leq
n-1$, there exist functions $f \in F_{\vec{n}}$, $g \in
G_{\vec{m}}$, such that $f$ and $g$ change sign exactly at each of
these points. This is a general property of Chebyshev spaces, see
\cite[Chapter 3.1, exercise E.11]{BE}.
\end{enumerate}

Now let $f \in F_{\vec{n}}$ be non-zero.
Then $f$ has at most $n-1$ real zeros of by (1). Let
$x_1, \ldots, x_m$ with $m \leq n-1$ be the zeros
of odd multiplicity (so that $f$ has a sign change
at these points). By (2) there is $g \in G_{\vec{m}}$
which also changes sign exactly at these points.
Then $fg$ has no sign change on $\mathbb R$, and
therefore $\int f(x) g(x) dx > 0$. Thus $f$ does not
belong to $G_{\vec{m}}^{\bot}$ and (\ref{Chebsys}) follows.
\end{proof}

By Lemmas \ref{Cheblemma} and \ref{bases} there exist biorthogonal bases
$\phi_1,\ldots,\phi_n$ of $F_{\vec{n}}$ and
$\psi_1,\ldots,\psi_n$ of $G_{\vec{n}}$.
Let $K_n$ be the projection kernel
\begin{align}
    K_n(x,y) = \sum_{j=1}^n \phi_j(x) \psi_j(y).
\end{align}
Then in the confluent case the probability density (\ref{karlin}) to
find the Brownian paths at time $t$ at the positions $x_1, \ldots, x_n$
can be written as
\[ p_{n,t}(x_1, \ldots, x_n) = \frac{1}{n!} \det( \phi_j(x_k))_{j,k=1}^n
    \det(\psi_j(x_k))_{j,k=1}^n = \frac{1}{n!} \det(K_n(x_j,x_k))_{j,k=1}^n.
    \]
Moreover, all correlation functions have determinantal form
with kernel $K_n$. That is, if
\[ r_m(x_1, \ldots, x_m) = \frac{n!}{(n-m)!}
    \underbrace{\int \cdots \int}_{n-m \mbox{ times}} p_{n,t}(x_1, \ldots, x_m, x_{m+1}, \ldots x_n) dx_{m+1} \cdots dx_n \]
denotes the $m$-point correlation function, then
\[ r_m(x_1, \ldots, x_m) = \det (K_n(x_j,x_k))_{j,k=1}^m \]
for every $m = 1, \ldots, n$.

There is a substantial literature on determinantal point processes
and non-intersecting random paths see e.g.\ the recent surveys
\cite{BKPV,Joh,OC,Sosh}
and references cited therein.

\medskip
Our Theorem \ref{theorem} relates the kernel $K_n$ to the
Riemann-Hilbert problem for multiple orthogonal polynomials
of mixed type. This opens up the possibility to analyze the
kernel in the large $n$ limit with the Deift/Zhou steepest descent
method for Riemann-Hilbert problems.

To obtain interesting limit behavior, one first modifies
the transition probability (\ref{transprob}) to
\begin{align} \label{transprob2}
    P_n(t,a,x) = \frac{\sqrt{n}}{\sqrt{2\pi t}} e^{-\frac{n}{2t}(x-a)^2}
\end{align}
so that the overall variance of the Brownian paths is reduced
with increasing $n$. With increasing $n$, the starting points $a_1, \ldots, a_p$ and the end points
$b_1, \ldots, b_q$ remain fixed while the corresponding
multiplicities $n_1, \ldots, n_p$ and $m_1, \ldots, m_q$
increase with $n$, such that the limits
\[ \lim_{n \to \infty} \frac{n_j}{n} \quad \mbox{ and }
    \lim_{n \to \infty} \frac{m_j}{n}
\]
exist and are positive.

\begin{figure}
    \begin{center}
         {\includegraphics[width=8cm,height=6cm]{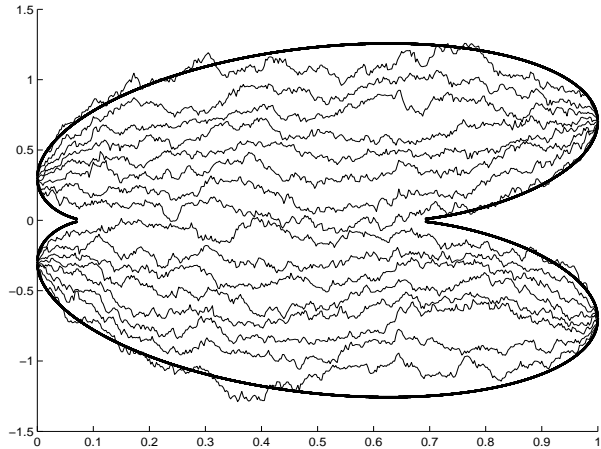}}
          \caption{\label{figuur4}
          Non intersecting Brownian motions which start and end at two different points}
    \end{center}
\end{figure}

Based on the experience with the case $p=1$ and $q=2$
that was developed in \cite{ABK3,BK2,BK4} we expect
the following to hold true.
In the limiting regime described above,
with probability one, the Brownian paths fill out a bounded
region as illustrated in Figure \ref{figuur4}.
Figure \ref{figuur4} shows a possible boundary curve
in the case of two starting points and two endpoints.
Here we see that two groups of paths start from two
different starting points and they come together
and merge at a certain critical time. Then they continue
as one group until at a second critical time they
split again into two groups that end at the two
different end points.
The boundary curve is smooth except for cusp singularities
that arise when two groups of paths come together or split.
We expect that this behavior is generic for general $p$ and $q$.

At any time $t \in (0,1)$ we further expect that the
correlation kernel $K_n$ has a scaling limit which
is equal to the usual scaling limits from random matrix theory.
That is, if we scale around a point $(t,x)$ lying strictly
inside the boundary curve then we expect the sine kernel
in the limit and for a usual point $(t,x)$ on the boundary
(not a cusp point) we expect the Airy kernel.

At the cusp singularities we expect that the kernel $K_n$
has the Pearcey kernel as a double scaling limit. This Pearcey kernel
arose first in the works of Br\'ezin and Hikami \cite{BH1,BH2} in the
context of Gaussian random matrices with external source.
In our notation this corresponds to $p=1$ and $q=2$.
A detailed treatment based on a double integral representation
of the kernel was made by Tracy and Widom \cite{TW}.
These authors also considered an extended Pearcey kernel and a
Pearcey process which involves the limiting joint distributions at
several scaled times near the critical times.
The Pearcey process also appears in the recent papers \cite{AvM,OR}.

For an extension of the above results to more general
values of $p$ and $q$ the Riemann-Hilbert problem that we gave in
this paper might be useful. Indeed, if the Deift/Zhou steepest descent
analysis can be made to work on this Riemann-Hilbert problem then the
scaling limits of the kernel can be derived.
We plan to report on this in a later publication.

\section*{Acknowledgements}

The authors are supported by FWO research projects G.0176.02 and
G.0455.04,  by K.U.Leuven research grant OT/04/24,
by INTAS Research Network NeCCA 03-51-6637,
by NATO Collaborative Linkage Grant PST.CLG.979738, by
grant BFM2001-3878-C02-02 of the Ministry of Science and
Technology of Spain, by the European Science Foundation Program
MISGAM and by the European Union through the FP6
Marie Curie RTN ENIGMA (Contract number MRTN-CT-2004-5652).

\end{document}